\documentclass[doublecol,linenumbers]{epl2} 
\usepackage{amsmath}
\usepackage{amssymb}
\usepackage{amsthm}
\usepackage{color}
\usepackage{indentfirst}

\title{Mean exit time and escape probability for the anomalous processes with the tempered power-law waiting times}
\shorttitle{Mean first exit time for the anomalous processes with the tempered waiting times} 

\author{ Weihua Deng, Xiaochao Wu, \and Wanli Wang}
\shortauthor{W.H. Deng, X.C. Wu, and W.L. Wang}

\institute{
   School of Mathematics and Statistics, Gansu Key Laboratory of Applied Mathematics and Complex
Systems, Lanzhou University, Lanzhou 730000,  P.R. China
}


\abstract{
The mean first exit (passage) time characterizes the average time of a stochastic process never leaving a fixed region in the state space, while the escape probability describes the likelihood of a transition from one region to another for a stochastic system driven by discontinuous (with jumps) L\'evy motion. This paper discusses the two deterministic quantities, mean first exit time and escape probability, for the anomalous processes having the tempered L\'{e}vy stable waiting times with the tempering index $\mu>0$ and the stability index $0<\alpha \le 1$; as for the distribution of jump lengths (in the CTRW framework) or the type of the noises driving the system (in the Langevin picture), two cases are considered, i.e., Gaussian white noise and non-Gaussian (tempered) $\beta$-stable ($0<\beta<2$) L\'{e}vy noise. Firstly, we derive the nonlocal elliptic partial differential equations (PDEs) governing the mean first exit time and escape probability. Based on the derived PDEs, it is observed that the mean first exit time depends strongly on the domain size and the values of $\alpha$, $\beta$ and $\mu$; when $\mu$ is close to zero, the mean first exit time tends to $\infty$. In particular, we also find an interesting result that the escape probability of a particle with (tempered) power-law jumping length distribution has no relation with the distribution of waiting times for the model considered in this paper.
For the solutions of the derived PDEs, the boundary layer phenomena are observed, which inspires the motivation for  developing the boundary layer theory for nonlocal PDEs.}

\begin{document}

\maketitle

\section{Introduction}
Anomalous diffusion phenomena are widely found in natural world; the subdiffusion includes, e.g.,  motion of lipids on membranes, solute transport in porous media, translocation of polymers; and the superdiffuion is observed in, e.g., turbulent flow, optical materials, motion of predators, human travel, etc \cite{Met:2000}. The types of diffusion are usually distinguished by the exponent of the evolution of the second order moment of a stochastic process $\mathbf{x}(t)$ with respect to the time $t$, i.e.,  $\langle \mathbf{x}^T(t) \mathbf{x}(t) \rangle \sim t^\gamma$; when $\gamma=1$, it is normal diffusion; $\gamma<1$ corresponds to subdiffusion and $\gamma>1$ superdiffusion.

 Anomalous diffusion is generally modelled by the continuous time random walks (CTRWs). Superdiffusion is usually associated with power law stepsize distribution and subdiffusion associated with power law waiting time distribution. The combination of power law stepsize and waiting time distributions may lead to superdiffusion, or subdiffusion, or even the normal diffusion, depending on the balance of the stepsize and waiting time distributions. For all these cases, the diffusion exponent does not change with time. Sometimes, the transition of the type of diffusions may occur over time, being detected in interplanetary solar-wind velocity and magnetic-field fluctuations \cite{Bruno:2004}, turbulent transport in magnetically confined plasmas \cite{Cartea:2007}, cage effect in a sheared granular \cite{Marty:2005}, diffusion of solar magnetics elements (being subdiffusive for times less than 20 minutes but normal for times larger than 25 minutes) \cite{Cadavid:1999}, motion of molecules diffusing in living cells \cite{Platani:2002}. The description of this kind of tempered anomalous dynamics can be realized by truncating the heavy tail of the power-law distribution \cite{Sokolov:2004} in the CTRW model. But here we use the exponential tempering, since it offers both mathematical and practical advantages, i.e., the tempered process is still an infinitely divisible L\'{e}vy process \cite{Meerschaert:2012, Sabzikara:2015}.

The first hitting time plays an important role in many fields, which is defined as the time when a certain condition is fulfilled by the random variable of interest for the first time. Generally, it is called first passage time when the random variable reaches a certain level for the first time, and first exit time when leaving a certain interval for the first time \cite{Martin:2011}. The example of first passage time naturally reaching our mind is the decision of an investor to buy or sell stock when its fluctuating prices reach a certain threshold. At the same time, a lot of research works reveal that anomalous process is a more effective model for the stock market. More examples of first passage time appear in chemical physics, e.g., certain chemical reaction occurs if a critical energy is reached through collisions or some other ways. In chemical physics, anomalous diffusion appears more often. All of these intrigue us to consider the first hitting time and the related issue: escape probability, for anomalous diffusion processes. In fact, for this topic there are already some relevant research works: Eliazar and  Klafter discussed the first passage leapovers based on the one-sided L\'{e}vy motion \cite{Eliazar2004On}; Bel and Barkai discussed the first passage time for unbiased and uniformly biased CTRW particles \cite{Bel2006Random}; Dybiec and his co-workers studied the mean first passage time on finite intervals for L\'{e}vy-Brownian motion \cite{Dybiec2006Levy}; Gajda  and Magdziarz studied the Kramers's escape problem for fractional Klein-Kramers equation \cite{Gajda2011Kramers}.
This paper is turning to derive the nonlocal partial differential equations (PDEs) governing the mean first exit time and escape probability for the tempered anomalous systems, based on the corresponding coupled Langevin equations \cite{Fogedby:1994}. According to the derived PDEs, more detailed (mathematical) analysis can be performed. Sometimes, it is very convenient.

In the past decades, most of the research works on the mean first exit time or escape probability, appearing in the mathematical, physical, chemical, and engineering literatures~\cite{Ben:1982,Bobrovsky:1982,Carmeli:1983, Gao:2014, Duan:2015, Day:1990,Friedman:1975,Gardiner:1985,Naeh:1990}, are for the uncoupled Langevin type dynamical system
\begin{equation}\label{Ito}
dX_{t}=F(X_{t})dt+\varepsilon\sigma(X_{t})dW_{t},
\end{equation}
where $W_t$ indicates the Gaussian or non-Gaussian $\beta$-stable type L\'{e}vy process, and $\varepsilon$ is a parameter that measures the strength of the noise. The first exit time from the spatial domain $D$ is defined as follows:
$$\tau(\mathbf{x}):=\inf\{t\geq0, X_{t}(\mathbf{x},\mathbf{y})\notin D\};$$
and the mean first exit time $u(\mathbf{x}):=\langle\tau(\mathbf{x})\rangle$.

\begin{figure}
  \centering
  \includegraphics[height=4.5cm,width=7.5cm]{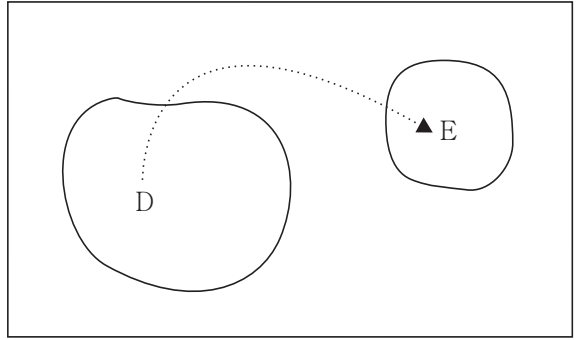}\\
  \caption{Sketch map of the escape probability.}\label{fig.n1}
\end{figure}

It is well known that if the noise matrix $\sigma(\cdot)$ has full rank throughout $D\cup\partial D$, the random trajectories of the solution of eq. (\ref{Ito}) hit $\partial D$ in finite time with probability 1 \cite{Friedman:1975}. Furthermore,
the random time $\tau$ to hit $\partial D$ for the first time has a finite first moment, and thus
\begin{equation}\label{finitefirstmom}
E\{\tau\}<\infty ~~~ {\rm and} ~~~   P_{r}(\tau<\infty)=1.
\end{equation}
The probability distribution of $\tau$, its moments, as well as the probability distribution of points on $\partial D$, where the random trajectories hit $\partial D$ for the first time, attract lots of interest in many applications \cite{Gardiner:1985}. If $W_t$ is a (tempered) non-Gaussian $\beta$-stable type L\'{e}vy process, we could also define another concept to quantify the exit phenomenon: escape probability, because of the discontinuity of the stochastic paths; the probability of a particle starting at a point $\mathbf{x}$, first escaping a domain $D$ and landing in a subset $E$ of $D^{c}$ (the complement of $D$), is called escape probability and denoted as $P_{E}(\mathbf{x})$; see fig. \ref{fig.n1}.

This paper considers the anomalous processes having the tempered L\'{e}vy stable waiting times with the tempering index $\mu>0$ and the stability index $0<\alpha \le 1$. We use the coupled Langevin equation to model the process, i.e., the stochastic trajectories of a $n$-dimensional continuous time random walk (CTRW) are expressed in terms of the coupled Langevin equations. The mean first exit time is discussed for two models, namely, the coupled Langevin equation driven by Gaussian white noise and the one by (tempered) non-Gaussian $\beta$-stable ($0<\beta<2$) L\'{e}vy noise. And the escape probability is analyzed for the second model. Firstly, we derive the nonlocal elliptic PDEs governing the deterministic quantities: the mean first exit time and the escape probability. The further investigations on the derived PDEs reveal that the mean first exit time depends strongly on the domain size and the values of $\alpha$, $\beta$ and $\mu$; in particular, when $\mu$ approaches to zero, the mean first exit time tends to $\infty$. Another interesting result is that the escape probability of a particle with power-law jumping length distribution is independent of the distribution of waiting times for the considered model in this paper.

All in all, this paper presents a direct method for the computation of the mean first exit time and the escape probability $P_{E}(x)$ for the tempered dynamical systems, in the limit of small Gaussian or non-Gaussian $\alpha$-stable type L\'{e}vy noises. And some interesting and general results are obtained. As for the derived PDEs, the boundary layer phenomena are observed, which inspires the motivation for developing the boundary layer theory for nonlocal PDEs.

%
%
%
%

\section{Gaussian white noise}\label{sect.1}
With the Gaussian white noise $\xi(t)$, the stochastic trajectory of a $n$-dimensional CTRW $Y(t)$ (subordinated Brownian process with external potential) is expressed in terms of the  coupled Langevin equation \cite{Fogedby:1994}:
\begin{equation}\label{Langevineq}
\begin{split}
&\dot{X}(s)=F(X(s))+\sqrt{2\varepsilon}\sigma(X(s))\xi(s),\\
&\dot{T}(s)=\eta(s),
\end{split}
\end{equation}
where $\eta(s)$ models the waiting times of the tempered anomalous diffusion process, being assumed to be independent from the $X$ process, i.e., the noises $\xi(s)$ and $\eta(s)$ are statistically independent \cite{Bouchaud:1990, BouchaudMF:1990}. $F(\mathbf{x})$ is a smooth vector field in a bounded domain $D$ in $\mathbf{R}^{n}$, whose boundary $\partial D$ has a smooth unit outer normal.
$F(\mathbf{x})$ and $\sigma(\mathbf{x})$ satisfy standard conditions \cite{Revuz1999} and the It\^{o} convention is adopted for the multiplicative term of eq. (\ref{Langevineq}).
 And $\sigma(\mathbf{x})$ is a $n\times k$ matrix of smooth noise coefficient, with $\varepsilon$ a parameter that measures the strength of the noise.
$\xi(s)$ is $k$-dimensional white Gaussian with $\langle\xi(s)\rangle=0$ and $\langle\xi^T(s_{1})\xi(s_{2})\rangle=\delta(s_{2}-s_{1})$ ($k\leq n$). The CTRW is then given by $Y(t)=X(S(t))$, where the process $S$ is defined as the inverse of $T$, or more precisely as the collection of first exit times:
\begin{equation}
S(t)=\displaystyle\inf_{s>0}\{s: \,T(s)>t\}.
\end{equation}
In the operational time $s$, the dynamics of $X$ is that of a normal diffusive process.
We take $\eta(s)$ as a tempered one-sided L\'{e}vy-stable noise with tempering index $\mu$ and stability index $0<\alpha<1$ interpolating between exponentially distributed $(\mu\rightarrow\infty)$ and power-law distributed $(\mu=0)$ waiting times \cite{Stanislavsky:2008,Gajda:2010}, obtained by the characteristic function of $T$: $\langle e^{-\lambda T(s)} \rangle=e^{-s ((\lambda+\mu)^\alpha-\mu^\alpha)}$.

 It can be noted that all the statistical properties of system eq. (\ref{Langevineq}) are determined by its transition probability density function (pdf) $$
 p(\mathbf{x},\mathbf{y},t)\mathbf{dy} \equiv Pr\{\mathbf{x(t)} \in \mathbf{y}+ \mathbf{dy} \,|\, \mathbf{x(0)}=\mathbf{x}  \},
  $$
which satisfies the
tempered fractional backward Kolmogorov equation \cite{Risken:1996,Magdziarz2009Langevin, Cairoli:2015} (for the forward version, see Appendix)
\begin{equation}\label{backwardKEq}
\frac{\partial}{\partial t}p(\mathbf{x},\mathbf{y},t)=\frac{\partial}{\partial t}\int_{0}^{t}K(t-t',\mu)L_{\mathbf{x}}^{*}p(\mathbf{x},\mathbf{y},t')dt',
\end{equation}
where the Laplace transform of the memory kernel is given by 
 $K(\lambda,\mu)=\frac{1}{(\lambda+\mu)^{\alpha}-\mu^{\alpha}}$ \cite{Cairoli:2015,Deng:2016};
the Laplacian operator
\begin{equation} \label{LevyOper1}
L_{\mathbf{x}}^{*}=\displaystyle\sum_{i=1}^{n}F^{i}(\mathbf{x})\frac{\partial}{\partial x^{i}}+\varepsilon\displaystyle\sum_{i,j=1}^{n}a^{ij}(\mathbf{x})\frac{\partial^{2}}{\partial x^{i}\partial x^{j}};
\end{equation}
%
%
 and the diffusion matrix $a(\mathbf{y})$ is defined by $a(\mathbf{y})=\frac{1}{2}\sigma(\mathbf{y})\sigma(\mathbf{y})^{T}$.
%
The solution of eq. ({\ref{backwardKEq}}), with the absorbing boundary condition
\begin{equation}\label{absboundcond}
p(\mathbf{x},\mathbf{y},t)|_{\mathbf{x}\in D,\,\mathbf{y}\in \partial D}=0
\end{equation}
and the initial condition
\begin{equation}\label{initcond}
p(\mathbf{x},\mathbf{y},0)=\delta(\mathbf{y}-\mathbf{x})
\end{equation}
has some special implications \cite{Naeh:1990}, as discussed below.

The condition (\ref{absboundcond}) implies that the particle is `absorbed' once it hits $\partial D$ for the first time, at time $t=\tau$; and it also signifies that
\begin{equation}\label{absboundcond2}
p(\mathbf{x},\mathbf{y},t)|_{\mathbf{x}\in\partial D,\,\mathbf{y}\in D}=0.
\end{equation}
The initial condition (\ref{initcond}) means that the trajectory starts at the point $\mathbf{x}$ in $D$ with probability 1.
Evidently, $\tau$, the first exit time to the boundary, is independent of the boundary behavior of the process. Thus, for instance, a process with absorption at the boundary hits $\partial D$ at the same time as a process without absorption, since up to time $\tau$ both evolve according to the same dynamics (\ref{Langevineq}). For a process with absorption at the boundary the probability of not exiting by time $t$ is identical to the probability of  finding it at time $t$ at some point $\mathbf{y}$ inside $D$, that is, the exit time distribution is given by
\begin{equation}\label{exittimedistr}
P_{r}\{\tau>t\,|\,\mathbf{x}(0)=\mathbf{x}\}=\int_{D}p(\mathbf{x},\mathbf{y},t)d\mathbf{y},
\end{equation}
where $p(\mathbf{x},\mathbf{y},t)$ is the solution of eq. (\ref{backwardKEq}) with (\ref{absboundcond}) and (\ref{initcond}) as the boundary and initial conditions. Now the mean first exit time of trajectories that start at $\mathbf{x}\in D$, is given by
\begin{equation}\label{aeq111}
\begin{split}
u(\mathbf{x}) & \equiv E(\tau\,|\,\mathbf{x}(0)=\mathbf{x}) \\
&=\int_{0}^{\infty}td_{t}[P_{r}(\tau<t\,|\,\mathbf{x}(0)=\mathbf{x})-1].
\end{split}
\end{equation}
Utilizing integration by parts and (\ref{finitefirstmom}) (see Appendix for more detailed explanations of the disappearance of boundary terms), there exists
\begin{equation}\label{MFPT1}
u(\mathbf{x})=\int_{0}^{\infty}P_{r}(\tau>t\,|\,\mathbf{x}(0)=\mathbf{x})dt;
\end{equation}
then combining (\ref{exittimedistr}) and (\ref{MFPT1}), we have
\begin{equation}\label{MFPT2}
u(\mathbf{x})=\int_{0}^{\infty}\int_{D}p(\mathbf{x},\mathbf{y},t)d\mathbf{y}dt.
\end{equation}
Defining
$$P(\mathbf{x},\mathbf{y},t)=\int_{0}^{t}p(\mathbf{x},\mathbf{y},t')dt' $$
and  taking Laplace transform $t\rightarrow \lambda$ lead to
$$P(\mathbf{x},\mathbf{y},\lambda)=\frac{p(\mathbf{x},\mathbf{y},\lambda)}{\lambda}.$$
Using the above formula and the final value theorem of Laplace transform ($\lim_{t\rightarrow\infty}f(t)=\lim_{\lambda\rightarrow0} \lambda f(\lambda)$), we get the stationary probability density
\begin{equation}\label{Pxy}
\begin{array}{lll}
\displaystyle P(\mathbf{x},\mathbf{y}) : = P(\mathbf{x},\mathbf{y},t=\infty)
&=& \displaystyle \lim_{\lambda\rightarrow 0}\lambda\cdot P(\mathbf{x},\mathbf{y},\lambda) \displaystyle
\\
&=& \displaystyle \lim_{\lambda\rightarrow 0}p(\mathbf{x},\mathbf{y},\lambda).
\end{array}
\end{equation}
The interpretation to the function $P(\mathbf{x},\mathbf{y})$ is the mean time that a trajectory $X(t)$ with $X(0)=\mathbf{x}$ spends at $\mathbf{y}$, prior to absorption at $\partial D$ \cite{Naeh:1990,Karlin:1981}.
Taking Laplace transform for eq. (\ref{backwardKEq}) with the initial condition (\ref{initcond}), we obtain
$$\lambda p(\mathbf{x},\mathbf{y},\lambda)=\delta(\mathbf{y}-\mathbf{x})+\lambda K(\lambda,\mu)L_{\mathbf{x}}^{*}p(\mathbf{x},\mathbf{y},\lambda). $$
Substituting $K(\lambda,\mu)$ into the above formula, and then letting $\lambda\rightarrow 0$ results in
\begin{equation}\label{Laplaceeq}
L_{\mathbf{x}}^{*}\displaystyle \lim_{\lambda\rightarrow 0}p(\mathbf{x},\mathbf{y},\lambda)=-\alpha\mu^{\alpha-1}\delta(\mathbf{y}-\mathbf{x}).
\end{equation}
Utilizing eq. (\ref{Pxy}) and (\ref{Laplaceeq}), we get
\begin{equation}\label{LaplaceeqforPxy}
L_{\mathbf{x}}^{*}P(\mathbf{x},\mathbf{y})=-\alpha\mu^{\alpha-1}\delta(\mathbf{y}-\mathbf{x});
\end{equation}
$P(\mathbf{x},\mathbf{y})$ also satisfies another forward type equation; see Appendix.

\begin{figure}
  \centering
  \includegraphics[height=5cm,width=8.5cm]{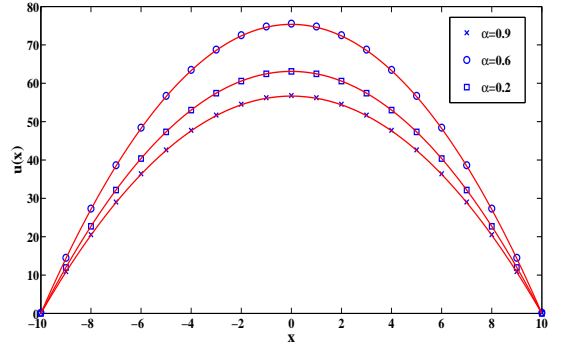}\\
  \caption{Behaviors of $u(x)$ generated by at least $5.5\times10^{5}$ trajectories for Gaussian jumping length with $\mu=0.1$, $n=1$, $r=10$. The (red) solid lines  are theoretical results (eq. (\ref{MFPTasy})). The symbols  (below with `$\times$' for $\alpha=0.9$, middle with `$\square$' for $\alpha=0.2$, above with `$\circ$' for $\alpha=0.6$) are for the simulation results. }\label{fig.1}
\end{figure}

Assuming that the order of the integral in eq. (\ref{MFPT2}) can be exchanged, then
\begin{equation}\label{MFPT3}
u(\mathbf{x})=\int_{D}\int_{0}^{\infty}p(\mathbf{x},\mathbf{y},t)dtd\mathbf{y}=\int_{D}P(\mathbf{x},\mathbf{y})d\mathbf{y}.
\end{equation}
Operating $L_{\mathbf{x}}^{*}$ on both sides of eq. (\ref{MFPT3}), we obtain
\begin{equation}
L_{\mathbf{x}}^{*}u(\mathbf{x})=\int_{D}L_{\mathbf{x}}^{*}P(\mathbf{x},\mathbf{y})d\mathbf{y}
=-\int_{D}\alpha\mu^{\alpha-1}\delta(\mathbf{y}-\mathbf{x})d\mathbf{y}.
\end{equation}
Therefore, the mean first exit time (\ref{MFPT1}) is the solution of the following equation
\begin{equation}\label{Mfpteq}
L_{\mathbf{x}}^{*}u(\mathbf{x})=-\alpha\mu^{\alpha-1} ~~~{\rm for}~~~\mathbf{x} \in D
\end{equation}
with the boundary condition
\begin{equation}\label{Mfptbc}
u(\mathbf{x})=0 ~~{\rm for}~~\mathbf{x} \in \partial D,
\end{equation}
because of eq. (\ref{absboundcond2}).
According to (\ref{LevyOper1}), eq. (\ref{Mfpteq}) can be rewritten as
\begin{equation}\label{n-dequation}
\displaystyle\sum_{i=1}^{n}F^{i}(\mathbf{x})\frac{\partial}{\partial x^{i}}u(\mathbf{x})+\varepsilon\displaystyle\sum_{i,j=1}^{n}a^{ij}(\mathbf{x})\frac{\partial^{2}}{\partial x^{i}\partial x^{j}}u(\mathbf{x})=-\alpha\mu^{\alpha-1}.
\end{equation}
Firstly, based on eq. (\ref{n-dequation}), we consider the stochastic dynamics of (\ref{Langevineq}) in one dimension. We take  $a(x)=1$, $F(x)=0$, $\varepsilon=1$ and $D=(-r,r)$ with $r>0$. Using the boundary condition eq. (\ref{Mfptbc}), by solving eq. (\ref{n-dequation}), we analytically obtain the mean first exit time given by
\begin{equation}\label{MFPTasy}
u(x)=\alpha\mu^{\alpha-1}\frac{r^2-x^2}{2},
\end{equation}
which is not a monotone function w.r.t. $\alpha$ if $\mu<1$, being confirmed by the simulation of the trajectories of the particles; see fig. \ref{fig.1}.

For $\alpha\in(0,1)$, the mean first exit time is decreasing with the increase of $\mu$, being easily understood; this is because the average waiting times of each step is becoming smaller. Moreover, if $\mu\rightarrow 0$, the mean first exit time $u(x)$ is infinite, which is an expected result because of the divergence of the average waiting times. Especially, when $\alpha=1$, $u(x)=\frac{r^2-x^2}{2}$.

Similar behaviors appear for the two dimensional case.
In fact, for the stochastic dynamics (\ref{Langevineq}) with $n=2$. We take $a^{ij}(\mathbf{x})=\left(
                                                                                                                \begin{array}{cc}
                                                                                                                   1& 0 \\
                                                                                                                   0& 1 \\
                                                                                                                \end{array}
                                                                                                              \right)
$, $F(\mathbf{x})=0$, $\varepsilon=1$ and a circular domain $D=\{\mathbf{x}: |\mathbf{x}|<r\}$ with $r>0$. Making use of the boundary condition eq. (\ref{Mfptbc}), the analytical solution (the mean first exit time) for eq. (\ref{n-dequation}) is
\begin{equation}\label{MFPTasy2d}
u(\mathbf{x})=\alpha\mu^{\alpha-1}\frac{r^2-|\mathbf{x}|^2}{4}.
\end{equation}
\section{Non-Gaussian L\'{e}vy noise}
Stochastic dynamics driven by non-Gaussian L\'{e}vy noises also have attracted much attention recently \cite{Applebaum:2009,Schertzer:2001,Kou2004Generalized,Dybiec2006Levy,Hofmann2003Mean}. We will specifically consider the following coupled Langevin equation:
\begin{equation}\label{levynoise}
\begin{split}
&\dot{X}(s)=F(X(s))+\varepsilon\dot{L}(s),\\
&\dot{T}(s)=\eta(s),
\end{split}
\end{equation}
where $F$ is a vector field (or drift); $L=(L(t), t\geq0)$ is a L\'{e}vy process defined in a probability space $(\Omega,\mathcal{F},\mathbb{P})$ with characteristics $(\mathbf{a},\mathbf{b},\nu)$, L\'{e}vy symbol $\psi$ and
pseudo-differential generator $A$ defined in eq. (\ref{pseudo-differential operator2}) \cite{Applebaum:2009}; $\varepsilon$ and $\eta(s)$ are the same as the ones used in  (\ref{Langevineq}); $L(s)$ and $\eta(s)$ are statistically independent. Without loss of generality, we may assume that the sample paths of $L$ are almost surely right continuous with left limits.

The L\'{e}vy-Khinchin formula says that any L\'{e}vy process has a specific
form for its characteristic function \cite{Applebaum:2009,Peszat:2007}, i.e., for all $t\geq 0$, $\mathbf{v}\in \mathbb{R}^{n}$,
$$\mathbb{E}(e^{i(\mathbf{v},L(t))})=e^{t\psi(\mathbf{v})},$$
where
\begin{equation}\label{psi(u)}
\begin{split}
& \psi(\mathbf{v})=i(\mathbf{a},\mathbf{v})-\frac{1}{2}(\mathbf{v},\mathbf{b}\mathbf{v})\\
&+\int_{\mathbb{R}^{n}\backslash\{0\}}\left[e^{i(\mathbf{v},\mathbf{y})}-1-i(\mathbf{v},\mathbf{y})\chi_{\{|\mathbf{y}|<1\}}\right]\nu(d\mathbf{y}),
\end{split}
\end{equation}
where $\chi_{I}$ is the indicator function of the set $I$, $\mathbf{a}\in \mathbb{R}^{n}$, $\mathbf{b}$ is a positive definite symmetric $n\times n$ matrix and $\nu$ is a sigma-finite  L\'{e}vy measure on $\mathbb{R}^{n}\backslash \{0\}$,
satisfying the property
\begin{equation}\label{borelmeasure}
\int_{\mathbb{R}^{n}\backslash \{0\}}(\mathbf{y}^{T}\mathbf{y} \wedge1)\nu(d \mathbf{y})<\infty,
\end{equation}
or equivalently
$$\int_{\mathbb{R}^{n}\backslash \{0\}}\frac{\mathbf{y}^{T}\mathbf{y}}{1+\mathbf{y}^{T}\mathbf{y}}\nu(d\mathbf{y})<\infty.$$
For each $f\in C_{c}^{\infty}(\mathbb{R}^{n})$, $\mathbf{x}\in\mathbb{R}^{n}$, 
\begin{equation}\label{pseudo-differential operator1}
(Af)(\mathbf{x})=(2\pi)^{-n/2}\int_{\mathbb{R}^{n}}e^{i(\mathbf{v},\mathbf{x})}\psi(\mathbf{v})\hat{f}(\mathbf{v})d\mathbf{v},
\end{equation}
where $\hat{f}(\mathbf{v})$ is the Fourier transform of $f(\mathbf{x})$.
Combining eq.~(\ref{psi(u)}) with (\ref{pseudo-differential operator1}), we have
\begin{equation}\label{pseudo-differential operator2}
\begin{split}
&(Af)(\mathbf{x})=F^i(\partial_if)(\mathbf{x})+a^{i}(\partial_{i}f)(\mathbf{x})+\frac{1}{2}b^{ij}(\partial_{i}\partial_{j}f)(\mathbf{x})+\\
&\int_{\mathbb{R}^{n}\backslash \{0\}}\left[f(\mathbf{x}+\mathbf{y})-f(\mathbf{x})-\mathbf{y}^{i}(\partial_{i}f)(\mathbf{x})\chi_{\{|\mathbf{y}|<1\}}\right]\nu(d\mathbf{y}),
\end{split}
\end{equation}
being the pseudo-differential operator for $X(z)$ of the system eq.~(\ref{levynoise}) with $\nu(d\mathbf{y})=\frac{\beta \Gamma(\frac{n+\beta}{2})}{2^{1-\beta} \pi^{n/2}\Gamma(1-\beta/2)} |\mathbf{y}|^{-\beta-n}d \mathbf{y}$,
where the Einstein summation convention \cite{Applebaum:2009} is used, and will also be used in the following.



Performing almost the same analysis for eq. (36) of \cite{Deng:2016} as done in Sec. II-B-2 of \cite{Deng:2016} and taking $\gamma$ and $q$ of the corresponding equation  equal to zero, we can obtain that the transition pdf $p(\mathbf{x},\mathbf{y},t)$ for system eq.~(\ref{levynoise}) satisfies the tempered fractional Kolmogorov equation \cite{Deng:2016,Applebaum:2009}
\begin{equation}\label{forwardKEq2}
\frac{\partial}{\partial t}p(\mathbf{x},\mathbf{y},t)=\frac{\partial}{\partial t}\int_{0}^{t}K(t-t',\mu)(Ap)(\mathbf{x},\mathbf{y},t')dt'.
\end{equation}
Similar to the discussion of last section, using the initial condition (\ref{initcond}) and absorbing boundary condition
\begin{equation}\label{absboundcondLevy}
p(\mathbf{x},\mathbf{y},t)|_{\mathbf{x}\in D,\,\mathbf{y}\in D^c}=0,
\end{equation}
which also implies that
\begin{equation}\label{absboundcondLevy2}
p(\mathbf{x},\mathbf{y},t)|_{\mathbf{x}\in D^c,\,\mathbf{y}\in D}=0,
\end{equation}
we have
\begin{equation}\label{Pxy2}
AP(\mathbf{x},\mathbf{y})=-\alpha\mu^{\alpha-1}\delta(\mathbf{y}-\mathbf{x});
\end{equation}
then the mean first exit time $u(\mathbf{x})$ for an orbit starting at $\mathbf{x}$ from a bounded domain $D$, satisfying the following differential-integral equation:
\begin{equation}\label{MFPTeq}
Au(\mathbf{x})=-\alpha\mu^{\alpha-1}~~~ {\rm for} ~~~\mathbf{x} \in D
\end{equation}
and
\begin{equation}\label{dirichletcond}
u(\mathbf{x})=0~~~{\rm for}~~~ \mathbf{x} \in  D^{c},
\end{equation}
where the generator $A$ is defined in eq.~(\ref{pseudo-differential operator2}), and $D^{c}=\mathbb{R}^{n}\backslash D$ is the complement set of $D$. Note that eq.~(\ref{dirichletcond}) is a nonlocal {\em Dirichlet} condition for the entire exterior domain $D^{c}$ or a volume constraint \cite{Du:2012}; using this type of boundary conditions, including (\ref{absboundcondLevy}), is due to the discontinuity of the trajectories of L\'{e}vy process.

\begin{figure}
  \centering
  \includegraphics[height=5cm,width=8.5cm]{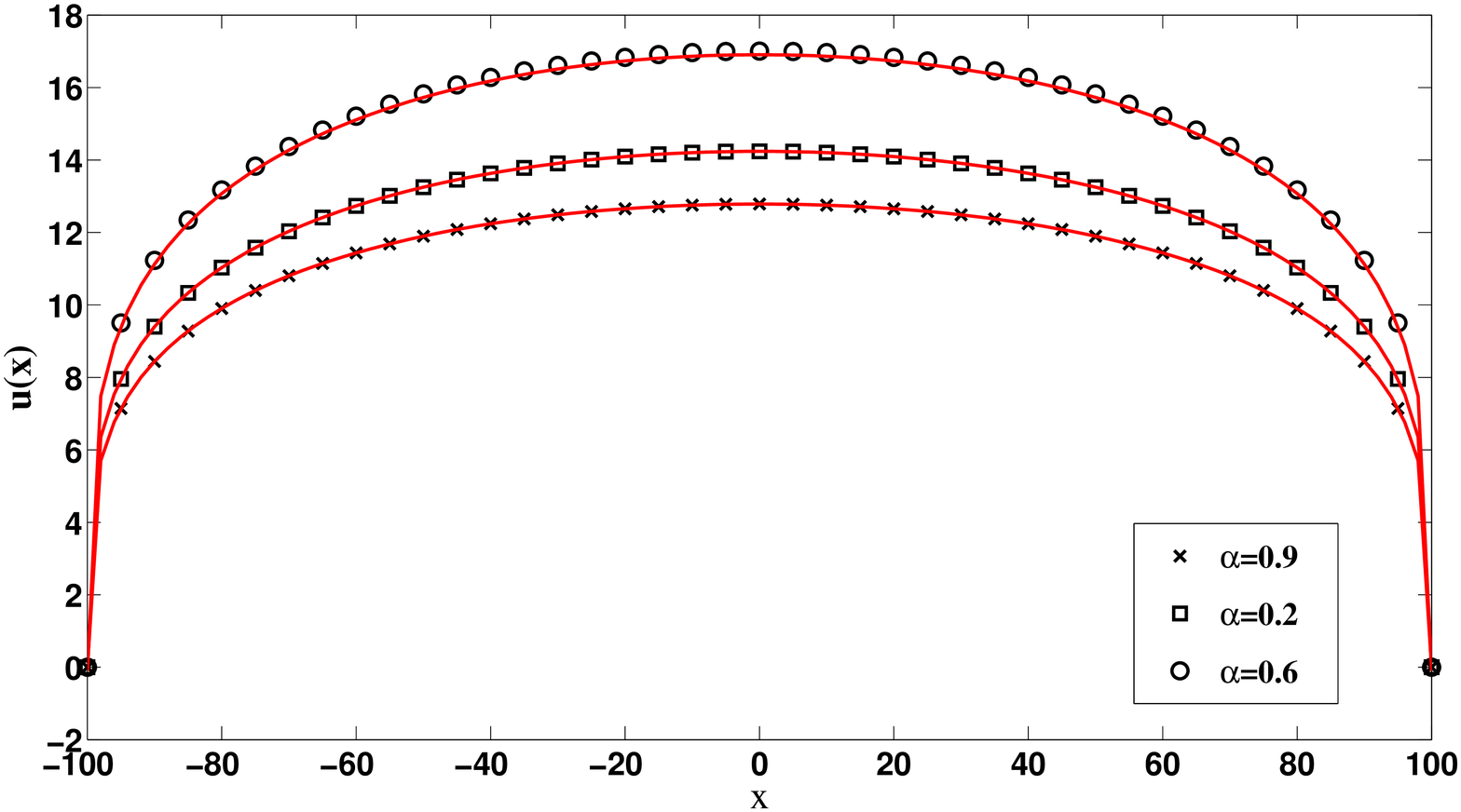}\\
  \caption{Behaviors of $u(x)$ generated by at least $1.5\times10^{6}$ trajectories with $\mu=0.1$, $\beta=0.5$, $n=1$, $r=100$. The (red) solid lines are theoretical results (eq. (\ref{MFPTasy2})).  The symbols  (below with `$\times$' for $\alpha=0.9$ , middle with `$\square$' for $\alpha=0.2$, above with `$\circ$' for $\alpha=0.6$) are for the simulation results. }\label{fig.2}
\end{figure}

\begin{figure}
  \centering
  \includegraphics[height=5cm,width=8.5cm]{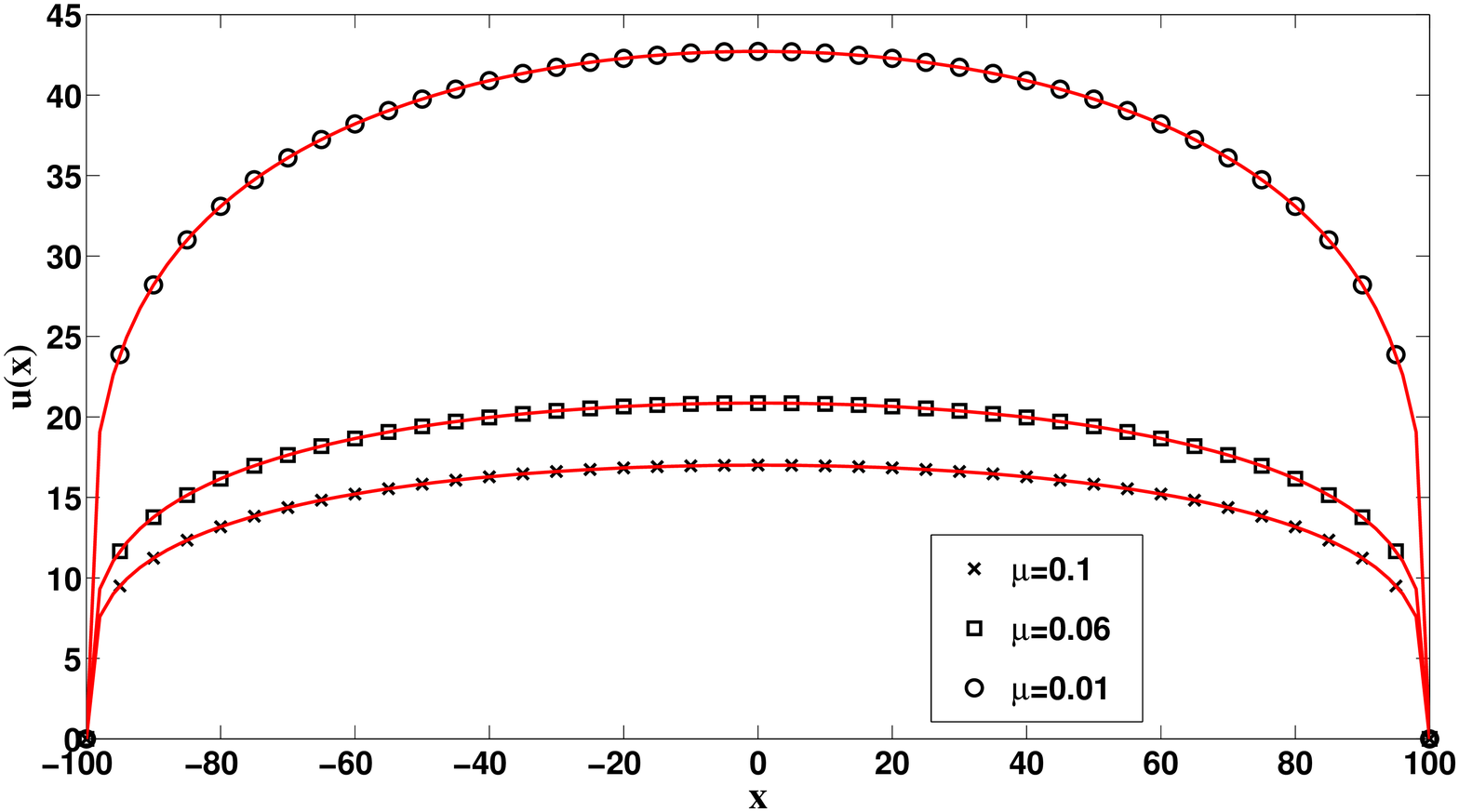}\\
  \caption{Behaviors of $u(x)$ generated by at least $1.5\times10^{6}$ trajectories with $\alpha=0.6$, $\beta=0.5$, $n=1$, $r=100$. The red solid lines  are theoretical results (eq. (\ref{MFPTasy2})). The symbols (below with `$\times$' for $\mu=0.1$, middle `$\square$' for $\mu=0.06$, above with `$\circ$' for $\mu=0.01$) are for the simulation results.}\label{fig.3}
\end{figure}

\begin{figure}
  \centering
  \includegraphics[height=4cm,width=8.5cm]{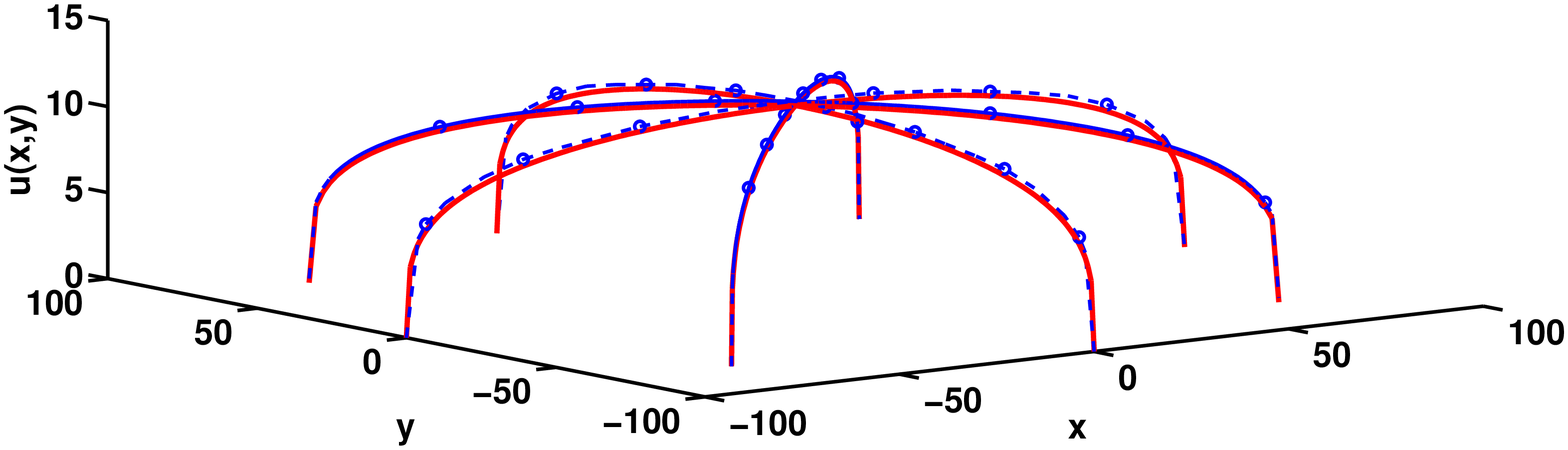}\\
  \caption{Behaviors of $u(x,y)$ generated by at least $2\times10^{6}$ trajectories with $\alpha=0.2$, $\mu=0.1$, $\beta=0.5$, $n=2$, $r=100$. The solid lines (red) are theoretical results (eq. (\ref{MFPTasy2})). The dotted (blue) lines with `$\circ$' are for the simulation results.}\label{fig.4}
\end{figure}

\begin{figure}
  \centering
  \includegraphics[height=4cm,width=8.5cm]{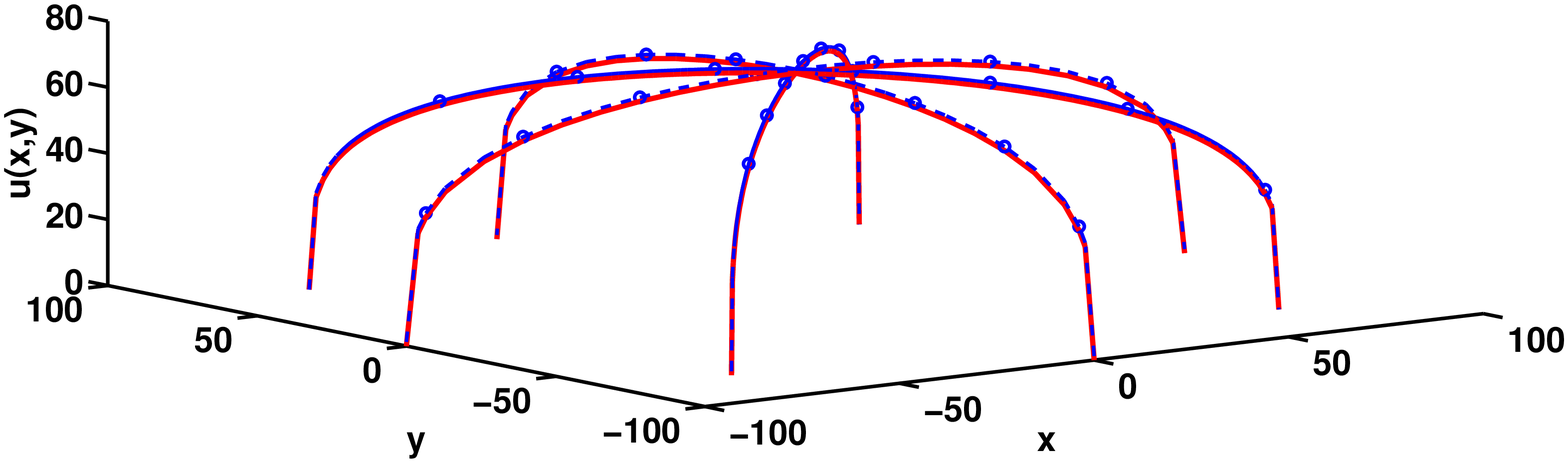}\\
  \caption{Behaviors of $u(x,y)$ generated by at least $2\times10^{6}$ trajectories with $\alpha=0.2$, $\mu=0.01$, $\beta=0.5$, $n=2$, $r=100$. The solid lines (red) are theoretical results (eq. (\ref{MFPTasy2})). The dotted (blue) lines with `$\circ$' are for the simulation results.}\label{fig.5}
\end{figure}

\begin{figure}
  \centering
  \includegraphics[height=4cm,width=8.5cm]{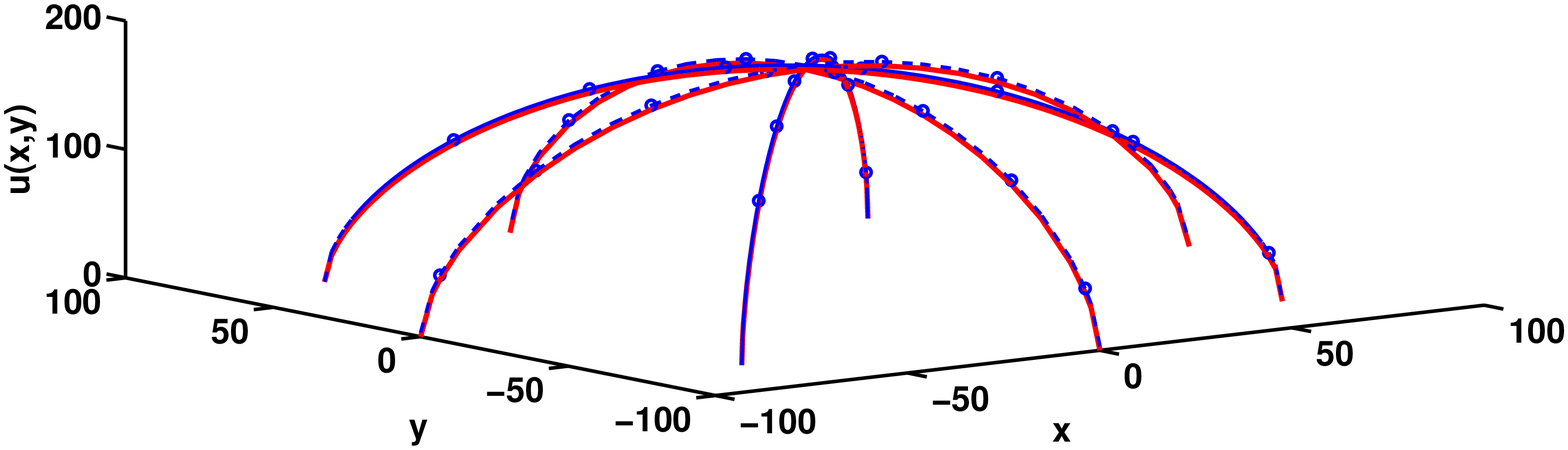}\\
  \caption{Behaviors of $u(x,y)$ generated by at least $2\times10^{6}$ trajectories with $\alpha=0.2$, $\mu=0.01$, $\beta=0.5$, $n=2$, $r=100$. The solid lines (red) are theoretical results (eq. (\ref{MFPTasy2})). The dotted (blue) lines with `$\circ$' are for the simulation results.}\label{fig.6}
\end{figure}

\begin{figure}
  \centering
  \includegraphics[height=4cm,width=8.5cm]{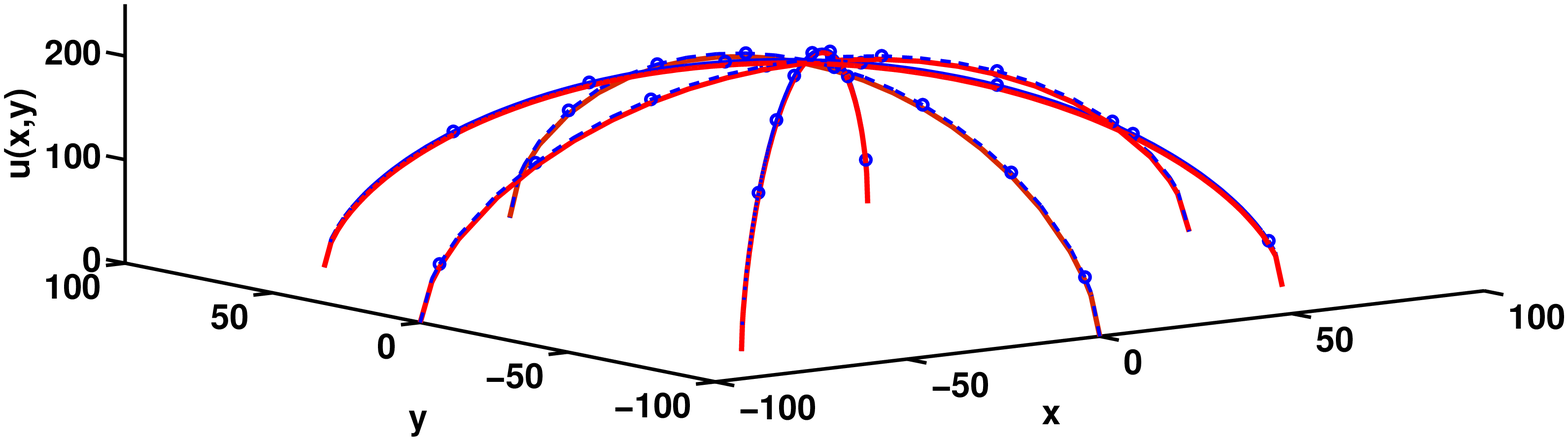}\\
  \caption{Behaviors of $u(x,y)$ generated by at least $2\times10^{6}$ trajectories with $\alpha=0.6$, $\mu=0.1$, $\beta=1.2$, $n=2$, $r=100$. The solid lines (red) are theoretical results (eq. (\ref{MFPTasy2})). The dotted (blue) lines with `$\circ$' are for the simulation results.}\label{fig.7}
\end{figure}

Considering a ball domain $D=\{\mathbf{x}:|\mathbf{x}|<r\}$, for $F(\mathbf{x})=0$ and $\varepsilon=1$, and a L\'{e}vy motion $L$ with the generating triplet $(0,0,\nu_{\beta})$, we have
\begin{equation}\label{MFPTasy2}
u(\mathbf{x})=\alpha\mu^{\alpha-1}\widetilde{u}(\mathbf{x}),
\end{equation}
where $\widetilde{u}(\mathbf{x})=\frac{\Gamma(n/2)(r^{2}-|\mathbf{x}|^{2})^{\beta/2}}{2^{\beta}\Gamma(1+\beta/2)\Gamma(n/2+\beta/2)}$ introduced in \cite{Getoor:1961}. It is obvious that the mean first exit time depends strongly on the domain size and the values of $\alpha$, $\beta$ and $\mu$; for the one dimensional case, see fig.~\ref{fig.2} and fig.~\ref{fig.3}; and for the two dimensional case, see fig. \ref{fig.4}, fig. \ref{fig.5}, fig. \ref{fig.6}, and fig. \ref{fig.7}. 
After making the comparison, especially the one among fig. \ref{fig.4}-\ref{fig.7}, the role of the parameters $\mu$, $\alpha$, and $\beta$ is more clearly shown, i.e., when $\mu$ tends to small, the mean exit time approaches to infinity and more stronger boundary layer phenomena are observed; when $\beta$ becomes smaller, the mean exit time tends to shorter but the boundary layer phenomena become stronger.
Moreover, for Gaussian jumping length $\beta=2$, taking $n=1$ or 2, then eq. (\ref{MFPTasy2}) expectedly corresponds to eq. (\ref{MFPTasy}) and (\ref{MFPTasy2d}), respectively.


If $\mu$ is finite and $\alpha=1$, eq.~(\ref{MFPTasy2}) reduces to $\widetilde{u}(\mathbf{x})$, i.e., truncation has no effect on normal diffusion, as predicted in \cite{Getoor:1961}. For $\alpha\in(0,1)$, the mean first exit time is increasing with the decrease of $\mu$. 
Moreover, if $\mu\rightarrow 0$, the mean first exit time $u(\mathbf{x})$ is infinite, because of the heavy-tailed distribution of the waiting times. For $n=1$ and $\alpha=\beta=1$, eq.~(\ref{MFPTasy2}) was obtained by Kac and Pollard \cite{Kac:1950}.

We then consider the escape probability of a particle whose motion is described by the coupled Langevin equation (\ref{levynoise}). The likelihood a particle, starting at a point $\mathbf{x}$, first escapes a domain $D$ and lands in a subset $E$ of $D^{c}$ is defined as escape probability. Basing on (\ref{Pxy2}), this escape probability, denoted by $P_{E}(\mathbf{x})$, taking advantage of initial condition (\ref{initcond}) and absorbing boundary condition (\ref{absboundcondLevy}),  satisfies
$$AP_{E}(\mathbf{x})=A\int_{E}P(\mathbf{x},\mathbf{y})d\mathbf{y}=\int_{E}-\alpha\mu^{\alpha-1}\delta(\mathbf{y}-\mathbf{x})d\mathbf{y}=0,$$
for all $\mathbf{x}\in D$. Hence, the escape probability $P_{E}(\mathbf{x})$ solves
\begin{equation}\label{Escapep}
\begin{split}
AP_{E}(\mathbf{x})=0,~~~~~~~~&\mathbf{x}\in D,\\
P_{E}(\mathbf{x})|_{\mathbf{x}\in E}=1,~~~~&P_{E}(\mathbf{x})|_{\mathbf{x}\in D^{c}\backslash E}=0,
\end{split}
\end{equation}
where $A$ is the generator also defined in (\ref{pseudo-differential operator2}). It can be easily observed that the escape probability $P_{E}(\mathbf{x})$ is independent of $\alpha$ and $\mu$, i.e., the waiting time distributions. And the solution to (\ref{Escapep}) with $D=(-1,1)$ and $E=[1,+\infty)$ is
$$
P_E(x)=2^{1-\alpha} \Gamma(\alpha)\left[ \Gamma\left( \frac{\alpha}{2}\right) \right]^{-2} \int_{-1}^x (1-u^2)^{\alpha/2-1}du,
$$
being the same as Corollary 1 of \cite{Blumenthal:1961}. Note that the concept of escape probability makes sense for all L\'{e}vy processes except Brownian motion (Gaussian noise), since almost all sample paths of Brownian motion are continuous in time in the common sense.


Considering the escape probability described by (\ref{Escapep}) in one dimension, in particular, we assume $D=(-r,r)$ with $r>0$ and $E=[r,\infty)$. Hence the conditions for the escape probability outside the domain are $P_{E}(x)=0$ for $x\in(-\infty,-r]$ and $P_{E}(x)=1$ for $x\in E$.
Meanwhile, we obtain the analytical result of the escape probability for the symmetric $\beta$-stable case $(F=0, a=0, b=0, \varepsilon=1)$ with $D=(-r,r)$ and $E=[r,\infty)$ \cite{Blumenthal:1961}:
\begin{equation}\label{PE}
P_{E}(x)=\frac{(2r)^{1-\beta}\Gamma(\beta)}{[\Gamma(\beta/2)]^{2}}\int_{-r}^{x}(r^2-y^2)^{\frac{\beta}{2}-1}dy,~~~x\in D.
\end{equation}
\begin{figure}
  \centering
  \includegraphics[height=5cm,width=8.5cm]{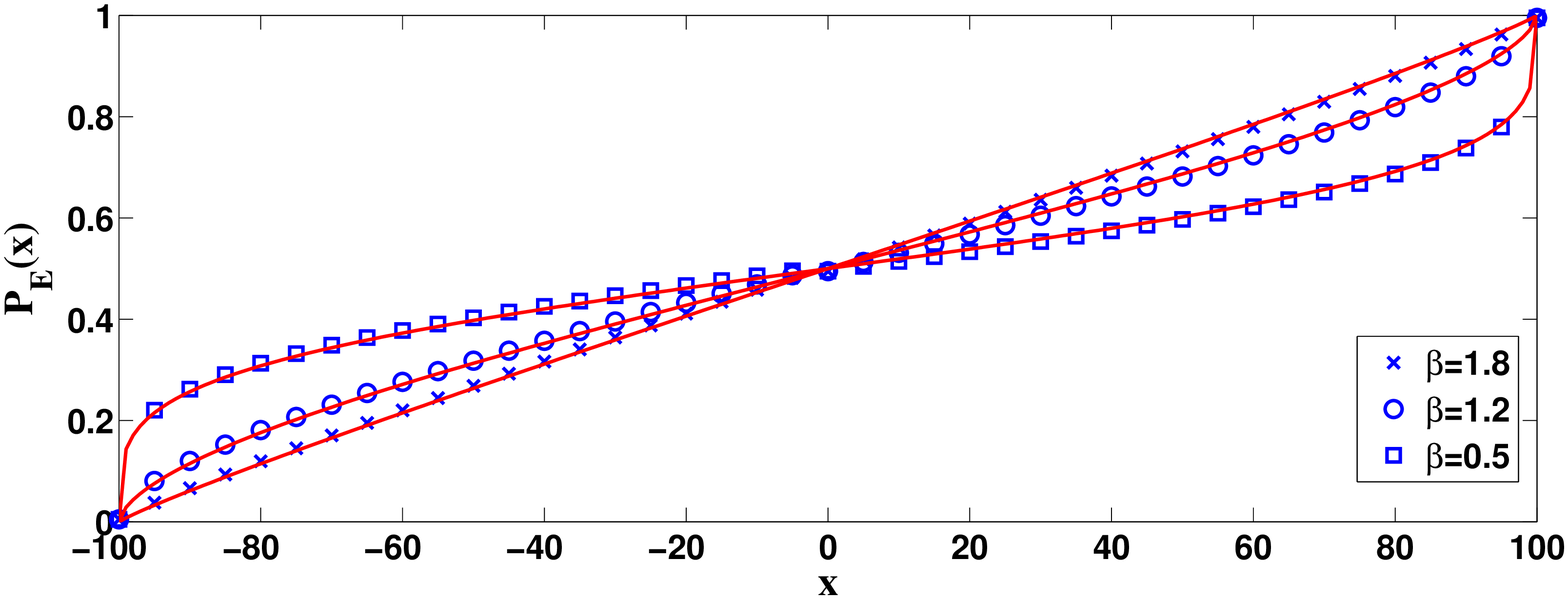}\\
  \caption{Behaviors of $P_{E}(x)$ generated by at least $10^{6}$ trajectories with $\alpha=0.6$, $\mu=0.1$, and $r=100$. The red solid lines  are theoretical results (eq. (\ref{PE})). The symbols denote random walk simulations using the trajectories of particles (left below right above with `$\times$' for $\beta=1.8$, middle with `$\circ$' for $\beta=1.2$, left above right below with `$\square$' for $\beta=0.5$.) }\label{fig.9}
\end{figure}

From the last equation, we can observe that due to the symmetry of the process and the domains, the escape probability $P_{E}(x)$ takes the value of one-half when the starting point is the position of symmetry $x=0$. The escape probability is symmetric with respect to the point in $(0,r)$, i.e.,
$$P_{E}(x)+P_{E}(-x)=1.$$
Because of the symmetry, we focus on positive starting points in the domain, i.e., $x>0$. It can be easily verified that for a fixed point, the probability for the process escaping to the right of the domain becomes smaller when the value of $\beta$ decreases. Moreover, this characteristic is independent of the domain size.
For the simulation results, see fig. \ref{fig.9}; it can be seen that at the vicinity of the boundary the curve becomes steeper when $\beta$ is smaller; boundary layer phenomena are observed.

\textbf{Remark:}
  If $L(s)$ is a tempered L\'evy process \cite{Meerschaert:2012}, the governing equations for the corresponding mean first exit time and escape probability are completely the same as eq. (\ref{MFPTeq}) and eq. (\ref{Escapep}), respectively, except replacing the operator $A$ there by
\begin{equation*}
\begin{split}
&(Af)(\mathbf{x})=F^{i}(\partial_{i}f)(\mathbf{x})+a^{i}(\partial_{i}f)(\mathbf{x})+\frac{1}{2}b^{ij}(\partial_{i}\partial_{j}f)(\mathbf{x})+\\
&\int_{\mathbb{R}^{n}\backslash \{0\}}\left[f(\mathbf{x}+\mathbf{y})-f(\mathbf{x})-\mathbf{y}^{i}(\partial_{i}f)(\mathbf{x})\chi_{\{|\mathbf{y}|<1\}}\right]\\
&C|\mathbf{y}|^{-\alpha-1}e^{-\mu |\mathbf{y}|}(d\mathbf{y}),
\end{split}
\end{equation*}
where $C$ is a normalized constant for the tempered L\'evy measure.

\section{Conclusion}
This paper focuses on the mean exit time and escape probability for the anomalous processes with the tempered power-law waiting times. Two models are considered: the one driven by Brownian motion (Gaussian noise) and the other one driven by (tempered) non-Gaussian $\beta$-stable ($0<\beta<2$) L\'{e}vy noise. The equations governing the mean exit times of the two models and escape probability of the latter model are derived. Based on the derived PDEs, two most striking results are obtained:  1. if the mean first exit time of the stochastic process $X(t)$ is $u(x)$, then it is $\alpha \mu^{\alpha-1} u(x)$ for the stochastic process $X(S(t))$, where $S(t)$ is the inverse of $T$ and $T(t)$ is the tempered L\'{e}vy stable process with the tempering index $\mu>0$ and the stability index $0<\alpha \le 1$; in fact, it holds for more general Markov processes (see the proof in Appendix); 2. the escape probability of a stochastic process is unrelated to the waiting time distribution of the stochastic process for the considered model Eq. (\ref{levynoise}). Other results include that the mean exit time tends to infinity if $\mu$ approaches zero and more detailed analyses on the dependencies of the mean exit time and escape probability on the parameters.

\acknowledgments
This work was supported by the National Natural Science Foundation of China under Grant No. 11671182. WD thanks Zhen-Qing Chen for the proof of the main result 1 described in the Conclusion section, being given in Appendix.
\section{Appendix}
\subsection{The tempered fractional forward Kolmogorov equation for the transition probability density function of eq. (\ref{Langevineq})}

\begin{equation}\label{forwardKEq}
\frac{\partial}{\partial t}p(\mathbf{x},\mathbf{y},t)=\frac{\partial}{\partial t}\int_{0}^{t}K(t-t',\mu)L_{\mathbf{y}}p(\mathbf{x},\mathbf{y},t')dt',
\end{equation}
where the Laplace transform of the memory kernel is given by 
 $K(\lambda,\mu)=\frac{1}{(\lambda+\mu)^{\alpha}-\mu^{\alpha}}$ \cite{Cairoli:2015,Deng:2016}, and  the Laplacian operator $L_{\mathbf{y}}=-\displaystyle\sum_{i=1}^{n}\frac{\partial}{\partial y^{i}}F^{i}(\mathbf{y})+\varepsilon\displaystyle\sum_{i,j=1}^{n}\frac{\partial^{2}}{\partial y^{i}\partial y^{j}}a^{ij}(\mathbf{y})$.

Similar to the derivation of eq. (\ref{LaplaceeqforPxy}), $P(\mathbf{x},\mathbf{y})$ defined in (\ref{Pxy}) also satisfies

\begin{equation}\label{LaplaceeqbackPxy}
L_{\mathbf{y}}P(\mathbf{x},\mathbf{y})=-\alpha\mu^{\alpha-1}\delta(\mathbf{y}-\mathbf{x}).
\end{equation}

\subsection{Derivation of eq. (\ref{MFPT1}) from eq. (\ref{aeq111})}
After performing integration by parts to eq. (\ref{aeq111}), we get eq. (\ref{MFPT1}) if the boundary terms disappear. Here we show it happens indeed. From eq. (\ref{aeq111}) we have
\begin{equation}\label{a2eq01}
\begin{split}
u(\mathbf{x})&=\int_{0}^{\infty}td_{t}[P_{r}(\tau<t\,|\,\mathbf{x}(0)=\mathbf{x})-1]\\
&=t[P_{r}(\tau<t\,|\,\mathbf{x}(0)=\mathbf{x})-1]|^{\infty}_{0}\\
&~~~~~+\int_0^{\infty}P_{r}(\tau>t\,|\,\mathbf{x}(0)=\mathbf{x})dt\\
&=\lim\limits_{t\rightarrow{\infty}}t[P_{r}(\tau<t\,|\,\mathbf{x}(0)=\mathbf{x})-1]\\
&~~~~~+\int_0^{\infty}P_{r}(\tau>t\,|\,\mathbf{x}(0)=\mathbf{x})dt.
\end{split}
\end{equation}
From eq. (\ref{finitefirstmom}), there exists
\begin{equation*}
  \int_0^{\infty}t\frac{\partial }{\partial t}[P_{r}(\tau<t\,|\,\mathbf{x}(0)=\mathbf{x})-1] dt<\infty,
\end{equation*}
which means 
\begin{equation}\label{a2eq02}
\lim\limits_{t\rightarrow{\infty}} \frac{t\frac{\partial }{\partial t}[P_{r}(\tau<t\,|\,\mathbf{x}(0)=\mathbf{x})-1]}{\frac{1}{t}}=0.
\end{equation}
According to L'H\^{o}pital's rule, there exists

\begin{equation}\label{a2eq03}
\begin{split}
&\lim\limits_{t\rightarrow{\infty}}t[P_{r}(\tau<t\,|\,\mathbf{x}(0)=\mathbf{x})-1]\\
&=\lim\limits_{t\rightarrow{\infty}}\frac{\frac{\partial }{\partial t}[P_{r}(\tau<t\,|\,\mathbf{x}(0)=\mathbf{x})-1]}{-\frac{1}{t^2}}=0.
\end{split}
\end{equation}
Combining eqs. (\ref{a2eq01}), (\ref{a2eq02}) and (\ref{a2eq03}), yields eq. (\ref{MFPT1}).

\subsection{Proof of the main result $1$ described in the Conclusion section}\label{app2}
Let $X_t$ be a general Markov process, and $Z_t$ be a strickly  increasing L\'{e}vy process, being independent of $X_t$. And $E_\tau$,  defined as $E_t=\inf\{\tau>0, Z_\tau>t\}$,  is an inverse subordinator. Furthermore we introduce the first exit time from the domain $D$: $\tau_D=\inf\{t>0, X_t\notin D\}$.  Letting $X_{t}^{*}=X_{E_t}$, we can obtain the following relations
\begin{equation*}\label{a3eq01}
\begin{split}
  \tau_D^{*}&=\inf\{t: X_t^{*}\notin D\}=\inf\{t: X_{E_t}^{*}\notin D\}\\
    & =\inf\{t: E_t>\tau_D\}=\inf\{t:  Z_{\tau_D}<t\}\\
    &= Z_{\tau_D}.
\end{split}
\end{equation*}
So, $E\{\tau_D^{*}\}=E\{Z_{\tau_D}\}=\int_0^\infty E\{Z_t\}P(\tau_D=t)dt$. Performing  Laplace transform  yields
\begin{equation*}
  E\{e^{-\lambda Z_t}\}=e^{-t\psi(\lambda)},
\end{equation*}
which implies that
\begin{equation*}
  E\{ Z_t\}=-\frac{\partial}{\partial \lambda}E\{ e^{-\lambda Z_t}\}|_{\lambda=0}=t\psi^{'}(0).
\end{equation*}
Hence, $E\{\tau_D^{*}\}=\int_0^{\infty}\phi^{'}(0)tP(\tau_D=t)dt=\phi^{'}(0)E\{\tau_D\}$. Then we get the final result
\begin{equation*}
  E\{\tau_D^{*}\}=\psi^{'}(0)E\{\tau_D\}.
\end{equation*}
Especially, letting $\psi(\lambda)=(\lambda+\mu)^\alpha-\mu^\alpha$, we have $E\{\tau_D^{*}\}=\alpha\mu^{\alpha-1}E\{\tau_D\}$, which agrees with the result we summarized in Conclusion section.


\begin{thebibliography}{0}

\bibitem{Met:2000}
  \Name{Metzler R. \and Klafter J.}
  \REVIEW{Phys. Rep.}{339}{2000}{1}.


\bibitem{Bruno:2004}
  \Name{Bruno R., Sorriso-Valvo L., Carbone V. \and Bavassano B.}
  \REVIEW{Europhys. Lett.}{66}{2004}{146}.

\bibitem{Cartea:2007}
  \Name{Cartea A. \and del-Castillo-Negrete D.}    
  \REVIEW{Phys. Rev. E}{76}{2007}{041105}.

\bibitem{Marty:2005}
  \Name{Marty G. \and Dauchot O.}
  \REVIEW{Phys. Rev. Lett.}{94}{2005}{015701}.

\bibitem{Cadavid:1999}
  \Name{Cadavid A. C., Lawrence J. K. \and Ruzmaikin A. A.}
  \REVIEW{Astrophys. J.}{521}{1999}{844}.

\bibitem{Platani:2002}
  \Name{Platani M., Goldberg I., Lamond A. I. \and Swedlow J. R.}
  \REVIEW{Nat. Cell Biol.}{4}{2002}{502}.


\bibitem{Sokolov:2004}
  \Name{Sokolov I. M., Chechkin A. V. \and Klafter J.}
  \REVIEW{Phys. A.}{336}{2004}{245}.

\bibitem{Sabzikara:2015}
  \Name{Sabzikara F., Meerschaerta M. M. \and Chen J. H.}
  \REVIEW{J. Comput. Phys.}{293}{2015}{14}.

\bibitem{Meerschaert:2012}
  \Name{Meerschaert M. M. \and Sikorskii A.}
  \Book{Stochastic Models for Fractional Calculus}
  \Publ{Walter de Gruyter \& Co, Berlin}
  \Year{2012}.



\bibitem{Martin:2011}
\Name{Martin E., Behn U. \and Germano G.}
\REVIEW{Phys. Rev. E.}{83}{2011}{051115}.

\bibitem{Eliazar2004On}
  \Name{Eliazar I. \and Klafter J.}
  \REVIEW{Physica A}{336}{2004}{219}.

\bibitem{Bel2006Random}
  \Name{Bel G. \and Barkai E.}
  \REVIEW{Phys. Rev. E}{73}{2006}{016125}.

\bibitem{Dybiec2006Levy}
  \Name{Dybiec B., Gudowska-Nowak E., \and H\"{a}nggi, P.}
  \REVIEW{Phys. Rev. E}{73}{2006}{046104}.



\bibitem{Gajda2011Kramers}
  \Name{Gajda, J. \and Magdziarz, M.}
  \REVIEW{Phys. Rev. E}{84}{2011}{021137}.


\bibitem{Fogedby:1994}
  \Name{Fogedby H. C.}
  \REVIEW{Phys. Rev. E.}{50}{1994}{1657}.



\bibitem{Ben:1982}
  \Name{Benjacob E., Bergman D. J., Matkowsky B. J. \and Schuss Z.}
  \REVIEW{Phys. Rev. A}{26}{1982}{2805}.









\bibitem{Bobrovsky:1982}
  \Name{Bobrovsky B. Z. \and Schuss Z.}
  \REVIEW{SIAM J. Appl. Math.}{42}{1982}{174}.

\bibitem{Carmeli:1983}
  \Name{Carmeli B. \and Nitzan A.}
  \REVIEW{Phys. Rev. Lett.}{51}{1983}{233}.





\bibitem{Day:1990}
  \Name{Day M. V.}
  \REVIEW{J. Math. Anal. Appl.}{147}{1990}{134}.


\bibitem{Gao:2014}
  \Name{Gao T., Duan J.,  Li X. \and Song R.}
  \REVIEW{SIAM J. Sci. Comp.}{36}{2014}{A887}.


\bibitem{Duan:2015}
  \Name{Duan J. Q.}
  \Book{An Introduction to Stochastic Dynamics}
  \Publ{Cambridge University Press, Cambridge, UK}
  \Year{2015}.




\bibitem{Friedman:1975}
  \Name{Friedman A.}
  \Book{Stochastic Differential Equations and Applications}
  \Publ{Academic Press, New York}
  \Year{1975}.

\bibitem{Gardiner:1985}
  \Name{Gardiner C. W.}
  \Book{Handbook of Stochastic Methods: For Physics, Chemistry and the Natural Sciences}
  \Publ{Springer-Verlag, New York}
  \Year{1985}.

\bibitem{Naeh:1990}
  \Name{Naeh T., Klosek  M. M., Matkowsky  B. J., \and Schuss Z.}
  \REVIEW{SIAM J. Appl. Math.}{50}{1990}{595}.

\bibitem{Bouchaud:1990}
  \Name{Bouchaud J. P. \and Georges  A.}
  \REVIEW{Phys. Rep.}{195}{1990}{127}.


\bibitem{BouchaudMF:1990}
  \Name{Shlesinger  M. F., Klafter  J., \and  Wong  Y. M. }
  \REVIEW{J. Stat. Phys.}{27}{1982}{499}.



\bibitem{Revuz1999}
  \Name{ Revuz  D. \and  Yor M.}
  \Book{Continuous Martingales and Brownian
Motion}
  \Publ{Springer Berlin Heidelberg, Berlin}
  \Year{1999}.







\bibitem{Stanislavsky:2008}
  \Name{Stanislavsky A., Weron K. \and Weron A.}
  \REVIEW{Phys. Rev. E.}{78}{2008}{051106}.

\bibitem{Gajda:2010}
  \Name{Gajda J. \and Magdziarz M.}
  \REVIEW{Phys. Rev. E.}{82}{2010}{011117}.

\bibitem{Risken:1996}
   \Name{Risken H. \and Frank T.}
   \Book{The Fokker-Planck equation: methods of solution and applications}
   \Publ{Springer-Verlag, New York}
   \Year{1996}.

\bibitem{Magdziarz2009Langevin}
   \Name{ Magdziarz, M. }
   \REVIEW{J. Stat. Phys.}{135}{2009}{763}.


\bibitem{Cairoli:2015}
  \Name{Cairoli A. \and Baule A.}
  \REVIEW{Phys. Rev. Lett.}{115}{2015}{110601}.


\bibitem{Deng:2016}
  \Name{Wu X. C., Deng W. H. \and Barkai E.}
  \REVIEW{Phys. Rev. E.}{93}{2016}{032151}.

\bibitem{Karlin:1981}
  \Name{Karlin S. \and Taylor H. M.}
  \Book{A Second Course in Stochastic Processes}
  \Publ{Academic Press, New York}
  \Year{1981}.



\bibitem{Applebaum:2009}
  \Name{Applebaum D.}
  \Book{L\'{e}vy Processes and Stochastic Calculus}
  \Publ{Cambridge University Press, Cambridge, UK}
  \Year{2009}.

\bibitem{Schertzer:2001}
  \Name{Schertzer D., Larcheveque M., Duan J., Yanovsky V. \and Lovejoy S.}
  \REVIEW{J. Math. Phys.}{42}{2001}{200}.

\bibitem{Kou2004Generalized}
  \Name{Kou, S. C.,  Xie, X. Sunney}
  \REVIEW{Phys. Rev. Lett.}{93}{2004}{180603}.



\bibitem{Hofmann2003Mean}
  \Name{Hofmann, Helmut, Ivanyuk, Fedor A}
  \REVIEW{Phys. Rev. Lett.}{90}{2003}{132701}.



\bibitem{Peszat:2007}
  \Name{Peszat S. \and Zabczyk J.}
  \Book{Stochastic Partial Differential Equations with L\'{e}vy Processes}
  \Publ{Cambridge University Press, Cambridge, UK}
  \Year{2007}.


\bibitem{Du:2012}
  \Name{Du Q., Gunzburger M., Lehoucq R. B. \and Zhou K.}
  \REVIEW{SIAM Rev.}{54}{2012}{667}.


\bibitem{Getoor:1961}
  \Name{Getoor R. K.}
  \REVIEW{Trans. Amer. Math. Soc.}{101}{1961}{75}.

\bibitem{Kac:1950}
  \Name{Kac M. \and Pollard H.}
  \REVIEW{Canad. J. Math.}{11}{1950}{375}.

\bibitem{Blumenthal:1961}
  \Name{Blumenthal R. M., Getoor R. K. \and Ray D. B.}
  \REVIEW{Trans. Amer. Math. Soc.}{99}{1961}{540}.



%

\end{thebibliography}
\end{document}